\documentclass{amsart}

\newtheorem{theorem}{Theorem}[section]
\newtheorem{lemma}[theorem]{Lemma}
\newtheorem{cor}[theorem]{Corollary}
\newtheorem{prop}[theorem]{Proposition}

\theoremstyle{definition}

\newtheorem{xca}[theorem]{Exercise}

\theoremstyle{remark}
\newtheorem{remark}[theorem]{Remark}

\numberwithin{equation}{section}

\def\vs#1{\vskip .#1 cm} 

\def\xrarrow{\xrightarrow} 

\def\<{\left<}
\def\>{\right>}

\newcommand{\field}[1]{\mathbb{#1}}
\newcommand{\ZZ}{\ensuremath{{\field{Z}}}}
\newcommand{\CC}{\ensuremath{{\field{C}}}}
\newcommand{\RR}{\ensuremath{{\field{R}}}}
\newcommand{\QQ}{\ensuremath{{\field{Q}}}}

\newcommand{\length}[1]{\ensuremath {\| #1 \|}}

\newcommand{\commentout}[1]{}

\def\ll{\lambda}

\newcommand{\cC}{\ensuremath{{\mathcal{C}}}}
\newcommand{\cD}{\ensuremath{{\mathcal{D}}}}

\newcommand{\cF}{\ensuremath{{\mathcal{F}}}}

\newcommand{\cH}{\ensuremath{{\mathcal{H}}}}

\newcommand{\cP}{\ensuremath{{\mathcal{P}}}}

\newcommand{\cT}{\ensuremath{{\mathcal{T}}}}

\def\a{\alpha}

\def\d{\partial}

\def\e{\epsilon}
\def\f{\phi}

\def\s{\sigma}
\def\Sig{\Sigma}
\def\t{\tau}

\def\z{\zeta}

\def\rel{{\rm\ rel\ }}

\def\ol{\overline}
\def\ul{\underline}
\def\wt{\widetilde}

\begin{document}

\title{Axioms for higher torsion invariants of smooth bundles}

\author{Kiyoshi Igusa}

\thanks{Supported by NSF Grants DMS 02-04386, DMS 03-09480.}




\begin{abstract} We explain the relationship between various characteristic classes for smooth manifold bundles known as ``higher torsion'' classes. We isolate two fundamental properties that these cohomology classes may or may not have:
additivity and transfer. We show that higher Franz-Reidemeister torsion and higher
Miller-Morita-Mumford classes satisfy these axioms. Conversely, any characteristic class of smooth bundles satisfying the two axioms must be a linear combination of these two examples. 

We also show how higher torsion invariants can be computed using only the axioms. Finally, we explain the conjectured formula of S. Goette relating higher analytic torsion classes and higher Franz-Reidemeister torsion.
\end{abstract}

\maketitle


\section*{Introduction}


Higher analogues of Reidemeister torsion and Ray-Singer analytic torsion were developed by J. Wagoner, J.R. Klein, the author, M. Bismut, J. Lott, W. Dwyer, M. Weiss, E.B. Williams, S. Goette and many others (\cite{[Wagoner:higher-torsion]}, \cite{[K:thesis]}, \cite{[IK1:Borel2]}, \cite{[Bismut-Lott95]}, \cite{[DWW]}, \cite{[BG2]}, \cite{[Goette01]}, \cite{[Goette03]}, \cite{[I:BookOne]}).

This paper develops higher torsion from an axiomatic viewpoint. There are three main objectives to this approach:
\begin{enumerate}
	\item Simplify the computation of these invariants.
	\item Isolate the key properties of higher torsion.
	\item Explain the theorems relating higher Franz-Reidemeister torsion, Miller-Morita-Mumford (tautological) classes and higher analytic torsion classes.
\end{enumerate}
The following two theorems are examples of results which will make much more sense from the axiomatic viewpoint.
\begin{theorem}[Hain,I,Penner]\label{thm01}
The higher Franz-Reidemeister torsion invariants for the Torelli group $\t_{2k}(T_g)\in H^{4k}(T_g;\RR)$ are proportional to the Miller-Morita-Mumford classes.
\end{theorem}
This had been conjectured by J.R. Klein \cite{[K:Torelli]}. The precise proportionality constant was computed in \cite{[I:BookOne]}. We will see that this theorem is an example of the uniqueness theorem for even higher torsion invariants.

The next theorem is Theorem 0.2 in \cite{[Goette01]}. See \cite{[BG2]},  \cite{[Goette01]},  \cite{[Goette03]} for more details.

\begin{theorem}[S. Goette]\label{thm02}
Suppose that $p:M\to B$ is a smooth bundle with closed oriented manifold fiber $X$, $\cF$ is a Hermitian coefficient system on $M$ and $H^\ast(X;\cF)$ admits a $\pi_1B$ invariant metric. Suppose further that there exists a fiberwise Morse function on $M$. Then the Chern normalizations (from \cite{[BG2]}) of the higher analytic torsion classes $\cT_{2k}(E)\in H^{4k}(B;\RR)$ and the higher Franz-Reidemeister torsion class are defined and agree up to a correction term which is a multiple of the transfer to $B$ of the Chern character of the vertical tangent bundle of $M$:
\[
	\cT_{2k}(E)=\t_{2k}^{ch}(E)+\z'(-2k)rk(\cF) tr_B^E(ch_{4k}(T^vM)).
\]
\end{theorem}

This theorem, together with the uniqueness theorem for odd higher torsion invariants proved below, suggests that nonequivariant higher analytic torsion classes are odd torsion invariants. Recently, Sebastian Goette has claimed that he can prove his theorem in general, i.e., without the existence of a fiberwise Morse function.

In this paper we define a \emph{higher torsion invariant} to be a characteristic class $\t\in H^{4k}(B;\RR)$ of ``unipotent'' smooth bundles $E\to B$ which satisfies two axioms. We show that each such invariant is the sum of even and odd parts
$
    \t=\t^++\t^-
$. The main theorem (Theorem \ref{thm:tau even and tau odd are unique up to scalar} below) is:

\begin{theorem}
Nontrivial even and odd torsion invariants $\t^+, \t^-$ exist in degree $4k$ for all $k>0$ and they are uniquely determined up to scalar multiples.
\end{theorem}

The uniqueness statement is simply a reflection of the fact that, given the axioms, the higher torsion invariant is easy to compute in many cases. We will carry out this computation for any unipotent smooth bundle pair $(E,\d_0E)\to B$ which admits a fiberwise Morse function $(E,\d_0E)\to (I,0)$ with the additional property that the critical points have distinct critical values.

To prove the existence we show that there are two linearly independent
higher torsion theories given by the higher Miller-Morita-Mumford
classes 
\[
	M_{2k}(E)=tr_B^E((2k)!ch_{4k}(T^vE))\in H^{4k}(B;\ZZ)
\]
and the higher Franz-Reidemeister (FR)
torsion invariants $\t_{2k}(E)\in H^{4k}(B;\RR)$. Using basic properties of higher
FR-torsion proved in \cite{[I:BookOne]} and
\cite{[I:ComplexTorsion]}, in particular the framing principle, it is easy to show that $\t_{2k}$ satisfies the
axioms. Basis properties of the transfer map imply that $M_{2k}$ also satisfies these axioms.

It is easy to see that $M_{2k}$ is an even higher torsion invariant, i.e., it is trivial when the fiber is a closed odd dimensional manifold. However, $\t_{2k}$ has both even and odd components. The uniqueness theorem implies that $\t_{2k}$ and $M_{2k}$ are proportional whenever the fiber is a closed oriented even dimensional manifold, for example an oriented surface. In order to determine the proportionality constant it suffices to compute one example.

Bismut and Lott \cite{[Bismut-Lott95]} showed that nonequivariant analytic torsion classes are trivial for bundles with closed even dimensional fibers. Thus, we believe that they are odd higher torsion invariants. S. Goette's theorem above says that, assuming the existence of a fiberwise Morse function, $\cT_{2k}$ is proportional to the odd part of $\t_{2k}$. And the proportionality constant is 1 if they are normalized in the same way.

An earlier version of this paper was entitled ``Axioms for higher torsion II.'' The present paper also incorporates relevant elements from my lecture notes
\cite{[I:AxiomsI]}. Thus this paper replaces both of these earlier works.

I am in great debt to John R. Klein and E. Bruce Williams for their help in
completing the final crucial steps in the proof of the main
theorem. I also benefitted greatly from conversations with
Sebastian Goette, Xiaonan Ma and Wojciech Dorabiala.


\section{Preliminaries}


We consider smooth fiber bundles
\[
	F\to E\xrarrow{p} B
\]
where $E$ and $B$ are compact smooth manifolds, $p$ is a smooth submersion and $F$ is a compact orientable manifold with or without boundary. In the boundary case there is a subbundle $\d F\to \d^vE\to B$ of $E$. We call $\d^vE$ the \emph{vertical boundary} of $E$. (The boundary of $E$ is the union of $\d^vE$ and $p^{-1}(\d B)$.)
We assume that $B$ is connected. We assume that the action of $\pi_1B$ on $F$ preserves some orientation of $F$.

We will assume that the bundle $E\to B$ is \emph{unipotent} in the sense that the rational homology of its fiber $F$ is unipotent as a $\pi_1 B$-module. In other words, $H_\ast(F;\QQ)$ has a filtration by $\pi_1B$ submodules so that the subquotients have trivial $\pi_1B$ actions. In particular, $\pi_1B$ does not permute the components of $F$. Note that unipotent $\pi_1B$ modules form a Serre category. In fact, it is the Serre category generated by the trivial modules.

\begin{prop}\label{prop:unipotent implies boundary unipotent}
If $E\to B$ is a unipotent bundle where $F$ has boundary $\d F$ then $H_\ast(\d F;\QQ)$ and $H_\ast(F,\d F;\QQ)$ are also unipotent $\pi_1B$ modules. In particular, the vertical boundary $\d^vE\to B$ is a unipotent bundle. 
\end{prop}

\begin{proof}
By Poincar\'{e} duality, $H_\ast(F;\QQ)\cong H^\ast(F;\d F;\QQ)$ is a unipotent $\pi_1B$ module. Its dual $H_\ast(F,\d F;\QQ)$ must also be unipotent. Since unipotent modules form a Serre category, the long exact homology sequence of $(F,\d F)$ implies that $H_\ast(\d F;\QQ)$ is also unipotent.
\end{proof}

Let $T^vE$ denote the \emph{vertical tangent bundle} of $E$. This is the bundle of all tangent vectors which lie in the kernel of $Tp:TE\to TB$. The \emph{Euler class} of the bundle
\[
    e(E)\in H^n(E,\d^vE;\ZZ)
\]
is simply the Euler class of $T^vE$. 

The \emph{transfer}
\begin{equation}\label{transfer}
    tr^E_B:H^\ast(E;\ZZ)\to H^\ast(B;\ZZ)
\end{equation}
is given by
\[
    tr^E_B(x)=p_\ast(x\cup e(E))
\]
where
\[
    p_\ast:H^{\ast+n}(E,\d^vE;\ZZ)\to H^\ast(B;\ZZ)
\]
is the push-down operator (given over $\RR$ by integrating along fibers).

If the orientation of the fiber $F$ is reversed, both $e(E)$
and $p_\ast$ change sign. Thus, the transfer is
independent of the choice of orientation of $F$. For the basic properties of the transfer see \cite{[Becker-Gottlieb:Adams-conj]}. The main property that we need is that, for closed fibers $F$,
\[
	tr_B^E=(-1)^n tr_B^E.
\]
So, rationally, $tr_B^E=0$ if $n=\dim F$ is odd.


\section{Axioms}


We define a \emph{higher torsion invariant} (in degree $4k>0$) to be a real characteristic class $\t(E)\in H^{4k}(B;\RR)$ for unipotent smooth bundles $E\to B$ with closed orientable fibers satisfying the \emph{additivity} and \emph{transfer} axioms described below.

When we say that $\t$ is a ``characteristic class'' we mean it is
a natural cohomology class. I.e.,
\[
    \t(f^\ast E)=f^\ast(\t(E))\in H^{4k}(B';\RR)
\]
if $f^\ast E$ is the pull-back of $E$ along $f:B'\to B$. Naturality implies that $\t$ is zero for trivial bundles: $\t(B\times F)=0$.

\subsection{Additivity} If $E=E_1\cup E_2$ where $E_1,E_2$ are unipotent bundles over $B$ with the same vertical boundary $E_1\cap E_2=\d^vE_1=\d^vE_2$ then the \emph{Additivity Axiom} says that
\begin{equation}\label{additivity axiom}
    \t(E)=\tfrac12\t(DE_1)+\tfrac12\t(DE_2)
\end{equation}
where $DE_i$ is the fiberwise double of $E_i$.

This wording of the Additivity Axiom comes from Ulrich Bunke.

\subsection{Transfer} Suppose that $p:E\to B$ is a unipotent bundle with closed fiber $F$ and $q:D\to E$ is an oriented $S^n$ bundle (associated to an $SO(n+1)$ bundle over $E$). Then the \emph{Transfer Axiom} says that the higher torsion invariants $\t_B(D)\in H^{4k}(B;\RR)$ and $\t_E(D)\in H^{4k}(E;\RR)$ are related by the formula:
\begin{equation}\label{transfer axiom}
    \t_B(D)=\chi(S^n)\t(E)+tr^E_B(\t_E(D))
\end{equation}
Note that $\chi(S^n)=2$ or $0$ depending on whether $n$ is even or odd respectively.


\subsection{Examples}


As stated in the introduction, two examples of higher torsion invariants are the higher Miller-Morita-Mumford classes $M_{2k}(E)$ and the higher Franz-Reidemeister torsion invariants $\t_{2k}(E)$. 

The Miller-Morita-Mumford classes, for closed fiber $F$, are given by
\[
	M_{2k}(E)=tr_B^E((2k)!ch_{4k}(T^vE))
\]
where $ch_{4k}(T^vE)=\frac12ch_{4k}(T^vE\otimes\CC)$. Although this is an integral cohomology class (for $k>0$) we consider it as a real characteristic class. This invariant is defined for any smooth bundle with closed orientable fiber $F$. If $n=\dim F$ is odd then twice the transfer map $tr_B^E$ is zero.

\begin{prop}\label{prop:M2k is even}
$M_{2k}(E)=0$ for closed odd dimensional fibers $F$.\qed
\end{prop}

\begin{theorem}\label{thm:M4k is a torsion invariant} $M_{2k}$ is a higher torsion invariant for every $k\geq1$.
\end{theorem}

The higher FR torsion invariants (\cite{[I:BookOne]}, \cite{[I:ComplexTorsion]})
\[
	\t_{2k}(E,\d_0E)\in H^{4k}(B;\RR)
\]
are defined for any \emph{relatively unipotent} bundle pair $(E,\d_0 E)\to B$. By this we mean that the vertical boundary $\d^v E$ is a union of two subbundles $\d_0 E,\d_1E$ with the same vertical boundary $\d_0E\cap\d_1E=\d^v\d_0E=\d^v\d_1E$ and that the rational homology of the fiber pair $(F,\d_0F)$ of $(E,\d_0E)$ is a unipotent $\pi_1B$ module.

\begin{theorem}\label{thm:FR is a torsion invariant}
The higher FR torsion invariants $\t_{2k}$ are higher torsion invariants for unipotent bundles with closed manifold fibers.
\end{theorem}

These two theorems will be proved later.

J-M. Bismut and J. Lott \cite{[Bismut-Lott95]} constructed even differential forms on $B$ called \emph{analytic torsion forms}. In some cases these are closed and topological (i.e., independent, up to exact forms, of the metric and horizontal distribution used to defined the form). For example, if $F$ is a closed oriented manifold and $\pi_1B$ acts trivially on $H_\ast(F;\QQ)$, then they obtain a (nonequivariant) \emph{analytic torsion class} 
\[
	\cT_{2k}^{BL}(E)\in H^{4k}(B;\RR).
\]
They showed that

\begin{prop}\label{prop:anaytic torsion is odd}
$\cT_{2k}^{BL}(E)=0$ for closed even dimensional fibers $F$.
\end{prop}

We also have the following theorems of X. Ma and U. Bunke.

\begin{theorem}[Ma \cite{[Ma97]}]\label{thm: Ma's theorem}
$\cT_{2k}$ satisfies the transfer axiom.
\end{theorem}

\begin{theorem}[Bunke \cite{[Bunke:spheres]}]
Let $E\to B$ be the $S^{2n-1}$ bundle associated to an $U(n)$ bundle $\xi$ over $B$. Then
\[
	\cT_{2k}^{Bunke}(E)=\frac{(4k+1)!}{2^{4k}(2k)!}\z(2k+1)ch_{4k}(\xi)
\]
\end{theorem}
\begin{remark} Bunke uses a different normalization of the analytic torsion. We should multiply by $(2\pi i)^{-2k}$ to get the Bismut-Lott normalization:
\[
	\cT_{2k}^{BL}(E)=(-1)^k(2\pi)^{-2k}\frac{(4k+1)!}{2^{4k}(2k)!}\z(2k+1)ch_{4k}(\xi)
\]
\end{remark}

As we noted in the introduction, S. Goette has extended the definition of the higher analytic torsion class to the case when $\pi_1B$ acts orthogonally on $H_\ast(F;\QQ)$, i.e., preserving some metric. However, we need to extend it to the unipotent case.


\section{Statement}


We give the statement of the main theorem. We begin with the following elementary observations.

\begin{lemma}
For each $k>0$ the set of all higher torsion invariants $\t$ of degree $4k$ is a vector space over $\RR$.
\end{lemma}

\begin{proof}
The axioms are homogeneous linear equations in $\t$.
\end{proof}

\begin{lemma}
If $\t$ is a higher torsion invariant then so is $(-1)^n\t$ where $n=\dim F$ considered as a function of the bundle $E\to B$.
\end{lemma}

\begin{proof}
We need to show that $(-1)^n\t$ satisfies the axioms:
\[
    (-1)^n\t(E)=(-1)^n\tfrac12\t(DE_1)+(-1)^n\tfrac12\t(DE_2)
\]
\[
    (-1)^{m+n} \t_B(S(\xi))=(-1)^n\chi(S^m)\t(E)+(-1)^m tr^E_B(\t_E(S(\xi))).
\]
The additivity axiom (the first equation) is the same as before. The transfer axiom (the second equation) is the same as before if both $n=\dim F$ and $m$ are even. If one or both are odd then the terms on the right with the wrong sign are zero since $\chi(S^m)=0$ for odd $m$ and $tr^E_B=0$ for odd $n$.
\end{proof}

These lemmas imply that higher torsion invariants can always be expressed as a sum of even and odd parts: $\t=\t^++\t^-$ where
\[
    \t^+=\t^+=\frac{\t+(-1)^n\t}2,\quad \t^-=\t^-=\frac{\t-(-1)^n\t}2.
\]
The even torsion invariant can only be nontrivial for even dimensional fibers and the odd torsion invariant can only be nontrivial for odd dimensional fibers. 

The main theorem of this paper is the following.

\begin{theorem}[Main Threorem]\label{thm:tau even and tau odd are unique up to scalar}
Nontrivial even and odd torsion invariants exist in degree $4k$ for all $k>0$ and they are unique up to a scalar factor.
\end{theorem}

\begin{cor}
Every even torsion invariant is a scalar multiple of $M_{2k}$ and every odd torsion invariant is a scalar multiple of the odd part $\t_{2k}^-$ of the higher FR-torsion $\t_{2k}$.
\end{cor}

The scalar multiples can be determined as follows. Let $\ll$ be the universal $U(1)=SO(2)$ bundle over $\CC P^\infty$. Let $S^1(\ll)\to \CC P^\infty$ be the circle bundle associated to $\ll$ and $S^2(\ll)\to \CC P^\infty$ the associated $S^2$ bundle (the fiberwise suspension of $S^1(\ll)$). Given any odd torsion theory $\t^-$ in degree $4k$ we have
\[
	\t^-(S^1(\ll))\in H^{4k}(\CC P^\infty;\RR)\cong\RR
\]
Therefore, $\t^-(S^1(\ll))=2s_1\,ch_{4k}(\ll)$ for some $s_1\in\RR$. The statement is that $\t^-$ is uniquely determined by the scalar $s_1$. 

Similarly, $\t^+(S^2(\ll))=2s_2\,ch_{4k}(\ll)$. And $\t^+$ is uniquely determined by the scalar $s_2$. For example, we have the following calculations which will be explained later.

\begin{prop}\label{prop:calc of M2k on S2(xi)}
$M_{2k}(S^2(\ll))=2(2k)!ch_{4k}(\ll)$.
\end{prop}

Thus, $s_2=(2k)!$ for the higher Miller-Morita-Mumford classes $M_{2k}$ (and $s_1=0$).

\begin{prop}\label{prop:calc of t2k on S1(xi)}
$\t_{2k}(S^n(\ll))=(-1)^{k+n}\z(2k+1)ch_{4k}(\ll)$ where $\z(s)=\sum_{m>0} 1/m^s$ is the Riemann zeta function.
\end{prop}

So, $s_n=\frac12(-1)^{k+n}\z(2k+1)$ for the higher FR torsion invariants $\t_{2k}$.

The uniqueness of even torsion gives us the following.

\begin{cor}[\cite{[I:ComplexTorsion]}]
If $E\to B$ is a unipotent bundle with closed even dimensional fibers then
\[
	\t_{2k}(E)=\frac{(-1)^k\z(2k+1)}{2(2k)!}M_{2k}(E).
\]
\end{cor}

Theorem \ref{thm01} in the introduction is a special case of this corollary.

The uniqueness of odd torsion can now be expressed as follows.

\begin{cor}
Any odd torsion invariant is a scalar multiple of the odd part of higher Franz-Reidemeister torsion which is given by
\[
	\t_{2k}^-=\t_{2k}-\frac{(-1)^k\z(2k+1)}{2(2k)!}M_{2k}.
\]
\end{cor}

Theorem \ref{thm02} in the introduction says that analytic torsion classes are odd torsion invariants on certain bundles. We expect that the same formula should hold in general.

The rest of this paper is devoted to the proof of the main theorem. We show that, given the values of the scalars $s_1,s_2$ above, any higher torsion invariant can be computed in sufficiently many cases to determine it completely.


\section{Extension to relative case}


In order to compute the higher torsion invariant $\t(E)$ we need to cut $E$ into simpler pieces and compute the relative torsion of each piece. To do this we need to extend $\t$ first to the case when $F$ has a boundary and then to the case of a unipotent bundle pair $(F,\d_0F)\to (E,\d_0E)\to B$.


\subsection{Higher torsion in the boundary case}


Suppose that $E\to B$ is a unipotent smooth bundle with fiber $F$ a compact orientable manifold with boundary. By Proposition \ref{prop:unipotent implies boundary unipotent} the vertical boundary $\d^vE$ is also unipotent. And it follows from the Mayer-Vietoris sequence that the fiberwise double $DE$ is also unipotent. A higher torsion invariant $\t$ can now be extended to the boundary case by the formula
\[
    \t(E):=\tfrac12\t(DE)+\tfrac12\t(\d^vE).
\]
We will show that this extension of $\t$ satisfies boundary analogues of the additivity and transfer axioms. We need the following lemmas.

\begin{lemma}\label{my wording of the additivity axiom} Suppose that $E_i$ are smooth unipotent bundles over $B$ with the same vertical boundary. Then
\[
    \t(E_1\cup E_2)+\t(E_3\cup E_4)=\t(E_1\cup E_3)+\t(E_2\cup E_4)
\]
\end{lemma}

\begin{proof}
Both sides are equal to $\frac12\sum \t(DE_i)$.
\end{proof}

\begin{lemma}\label{lem:torsion of dv E times D2}
$\t(\d^vE)=\t(\d^v(E\times D^2))$ assuming $E$ is unipotent.
\end{lemma}

\begin{proof}
Since $\d^v(E\times D^2)=\d^vE\times D^2\cup E\times S^1$ we have by the additivity axiom that
\[
    \t(\d^v(E\times D^2))=\tfrac12\t(\d^vE\times S^2)+\tfrac12\t(DE\times S^1).
\]
But $\t(\d^vE\times S^2)=2\t(\d^vE)$ and $\t(DE\times S^1)=0$ by the transfer axiom.
\end{proof}

\begin{lemma}[additivity of transfer]\label{lem:additivity of transfer}
 If $E=E_1\cup E_2$ is a union of two smooth bundles along their common vertical boundary $\d^vE_1=\d^vE_2=E_1\cap E_2$ then
    \[
        tr^E_B(x)=tr^{E_1}_B(x|E_1)+tr^{E_2}_B(x|E_2)-tr^{\d_0E_1}_B(x|\d^vE_1)
    \]
for all $x\in H^\ast(E;\RR)$.
\end{lemma}




\begin{prop}[additivity for boundary case]\label{additivity in boundary case} If $(E_1,\d_0), (E_2,\d_0)$ are unipotent bundle pairs over $B$ with $E_1\cap E_2=\d_0E_1=\d_0E_2$ then
    \[
        \t(E_1\cup E_2)=\t(E_1)+\t(E_2)-\t(E_1\cap E_2).
    \]
\end{prop}

\begin{proof} We expand each term using the defining equation:
\begin{align*}
   \t(E_i)&:=\tfrac12\t(DE_i)+\tfrac12\t(\d^v E_i)\\
    \t(E_1\cup E_2)&:=\tfrac12\t(D(E_1\cup E_2))+\tfrac12 \t(\d^v( E_1\cup
    E_2))\\
    \t(E_1\cap E_2)&:=\tfrac12\t(\d^v(E_1\cap E_2))+\tfrac12\t(D(E_1\cap
    E_2)).
\end{align*}
The order of the terms in the last equation is reversed so that it
matches the following two examples of Lemma \ref{my wording of the additivity axiom}:
\[
    \t(\d^vE_1)+\t(\d^vE_2)=
    \t(\d^v( E_1\cup E_2))+\t(D(E_1\cap E_2))
\]
\[
    \t(DE_1)+\t(DE_2)=\t(D(E_1\cup E_2))+\t(D(I\times (E_1\cap
    E_2))).
\]
However, the last term is
\[
    \t(D(I\times (E_1\cap
    E_2)))=\t(\d^v(D^2\times (E_1\cap
    E_2)))=\t(\d^v(E_1\cap E_2))
\]
by Lemma \ref{lem:torsion of dv E times D2}. The proposition follows.
\end{proof}




\begin{prop}[transfer for boundary case]\label{transfer in relative case} If $X\to D\to E$ is an oriented linear disk or sphere bundle then
    \[
        \t_B(D)=\chi(X)\t(E)+tr^E_B(\t_E(D)).
    \]
\end{prop}

\begin{proof}
We consider first the case when $D=D(\xi)$ is an oriented linear $D^n$-bundle and $F$ is closed. I.e., we will show: $\t_B(D(\xi))=\t(E)+tr^E_B(\t_E(D(\xi)))$. This is just half of the sum of the following two examples of the original transfer axiom.
\[
    \t_B(S^n(\xi))=\chi(S^n)\t(E)+tr^E_B(\t_E(S^n(\xi)))
\]
\[
    \t_B(S^{n-1}(\xi))=\chi(S^{n-1})\t(E)+tr^E_B(\t_E(S^{n-1}(\xi)))
\]
The transfer axiom takes care of the case when $D$ is a sphere bundle and $F$ is closed. The remaining case when $\d F$ is nonempty is given by the following lemma.
\end{proof}

\begin{lemma}
With the fiber $X$ of $D\to E$ fixed, the transfer formula for $F$ closed implies the transfer formula for $F$ with boundary.
\end{lemma}

\begin{proof}
Write $DE=E\cup E'$ as the union of two copies of $E$ along its vertical boundary. Let $D,D'$ be two copies of $D$ with $D\cap D'=q^{-1}(\d^vE)$. Then the transfer formula $\t_B(D)=\chi(X)\t(E)+tr^E_B(\t_E(D))$ is half the sum of the following two transfer formulas with closed fibers $DF,\d F$ respectively.
\begin{align*}
    \t_B(D\cup D')&=\chi(X)\t(E\cup E')
    +tr^{E\cup E'}_B(\t_{E\cup E'}(D\cup D'))\\
    \t_B(D\cap D')&=\chi(X)\t(E\cap E')
    +tr^{E\cap E'}_B(\t_{E\cap E'}(D\cap D'))
\end{align*}
The additivity of transfer (Lemma \ref{lem:additivity of transfer}) is used here.
\end{proof}


\subsection{Relative torsion} 


Suppose that $(F,\d_0F)\to (E,\d_0E)\to B$ is a \emph{unipotent smooth bundle pair}. By this we mean that the vertical boundary $\d^vE$ is the union of two subbundles $\d^vE=\d_0E\cup \d_1E$ which meet along their common vertical boundary: $\d_0E\cap \d_1E=\d^v\d_0E=\d^v\d_1E$ and that both $E$ and $\d_0E$ are unipotent. This implies the weaker condition that the pair $(E,\d_0E)$ is relatively unipotent. We use the abbreviation $(E,\d_0)$ for $(E,\d_0E)$.

Suppose that $\t$ is a higher torsion invariant which has been extended to the boundary case as above. Then for any unipotent smooth bundle pair $(E,\d_0)\to B$ we define the \emph{relative torsion} by
\[
    \t(E,\d_0):=\t(E)-\t(\d_0E).
\]




\begin{prop}[additivity in the relative case]\label{additivity in the relative case} Suppose that $E\to B$ is a smooth bundle which can be written as a union of two subbundles $E=E_1\cup E_2$ which meet along a subbundle of their respective vertical boundaries: $E_1\cap E_2=\d_0E_2\subseteq \d^vE_1$. Let $\d^vE_1=\d_0E\cup \d_1E$ be a decomposition $\d^vE_1$ so that $\d_0E_2\subseteq \d_1E_1$ and  $(E_i,\d_0)\to B, i=1,2$ are unipotent smooth bundle pairs. Then $(E,\d_0E_1)\to B$ is unipotent and
\[
    \t(E_1\cup E_2,\d_0E_1)=\t(E_1,\d_0)+\t(E_2,\d_0).
\]
\end{prop}

\begin{proof}
Both sides are equal to $\t(E_1)+\t(E_2)-\t(E_1\cap E_2)-\t(\d_0E_1)$.
\end{proof}

Here is another variation of the additivity axiom which is also trivial to prove.

\begin{prop}[horizontal additivity]\label{horizontal additivity} Suppose that $(E,\d_0)\to B$ is a union of two unipotent bundle pairs $(E_i,\d_0E_i)$ in the sense that $E=E_1\cup E_2$ and $\d_0E=\d_0E_1\cup \d_0E_2$ with $E_1\cap E_2\subseteq \d_1E_1\cap \d_1E_2$. Let $X=E_1\cap E_2$ and $\d_0X=X\cap \d_0E$ and suppose $(X,\d_0)$ is a unipotent bundle pair. Then $(E,\d_0)$ is unipotent and
\[
    \t(E,\d_0)=\t(E_1,\d_0)+\t(E_2,\d_0)-\t(X,\d_0).\qed
\]
\end{prop}

To state the transfer axiom in the relative case we need the \emph{relative transfer}:
\[
    tr^{(E,\d_0)}_B:H^\ast(E;\ZZ)\to H^\ast(B;\ZZ)
\]
given by
\begin{equation}\label{relative transfer}
    tr^{(E,\d_0)}_B(x)=p_\ast(x\cup e(E,\d_0))
\end{equation}
where $p_\ast$ is the push-down operator (as before) and
\[
    e(E,\d_0)\in H^n(E,\d^vE;\ZZ)
\]
is the \emph{relative Euler class} given by pulling back the Thom class of the vertical tangent bundle $T^vE$ along any vertical tangent vector field which is nonzero along the vertical boundary $\d^vE$ and which points inward along $\d_0E$, outward along $\d_1E$ and is tangent to $\d^vE$ pointing from $\d_0E$ to $\d_1E$ along $\d_0E\cap \d_1E$.

As in the absolute case, the relative transfer is independent of the choice of orientation of the fiber. The relative transfer also satisfies the following two equations for any $x\in H^\ast(E;\ZZ)$.
\begin{equation}\label{relative transfer equation}
	tr^{(E,\d_0)}_B(x)=tr^E_B(x)-tr^{\d_0E}_B(x|\d_0E)
\end{equation}
\begin{equation}\label{duality for relative transfer}
	 tr^{(E,\d_1)}_B(x)=(-1)^n tr^{(E,\d_0)}_{B}(x)
\end{equation}
And, finally, $tr_B^{(E,\d_0)}\circ p^\ast:H^\ast(B)\to H^\ast(B)$ is multiplication by the \emph{relative Euler characteristic} of the fiber pair $(F,\d_0F)$ given by
\[
	\chi(F,\d_0):=\chi(F)-\chi(\d_0F).
\]




\begin{prop}[transfer in the relative case]\label{transfer in the relative case}
Let $(F,\d_0)\to(E,\d_0)\xrarrow{p} B$ and $(X,\d_0)\to(D,\d_0)\xrarrow{q} E$ be unipotent smooth bundle pairs so that the second is an oriented linear $S^n$ or $D^n$ bundle with $\d_0X=S^{n-1},D^{n-1}$ or $\emptyset$. Then
\[
    \t_B(D,\d_0D\cup q^{-1}\d_0E)
    =\chi(X,\d_0)\t(E,\d_0)+tr^{(E,\d_0)}_B(\t_E(D,\d_0)).
\]
\end{prop}

\begin{proof} We already did the case when both $\d_0F$ and $\d_0X$ are empty.
The case when $\d_0X$ is empty follows easily from the following formula which holds by definition.
\[
    \t_B(D,q^{-1}\d_0E)=\t_B(D)-\t_B(q^{-1}\d_0E)
\]
The general case follows from the following two examples of the $\d_0X=\emptyset$ case.
\begin{align*}
    \t_B(\d_0D,\d_0D\cap q^{-1}\d_0E)
    &=\chi(\d_0X)\t(E,\d_0)+tr^{(E,\d_0)}_B(\t_E(\d_0D))\\
    \t_B(D,q^{-1}\d_0E)
    &=\chi(X)\t(E,\d_0)+tr^{(E,\d_0)}_B(\t_E(D))
\end{align*}
Take the second formula minus the first to prove the proposition.
\end{proof}

A useful special case is the case when $D\to E$ is the pull-back of a linear bundle over $B$. In this case $D$ is the fiber product of two bundles.




\begin{cor}[product formula]\label{cor:product formula}
Suppose that $(F,\d_0)\to (E,\d_0)\to B$ and $(X,\d_0)\to (E',\d_0)\to B$ are unipotent smooth bundle pairs so that the second is an oriented linear $S^n$ or $D^n$ bundle with $\d_0X=S^{n-1},D^{n-1}$ or $\emptyset$. Let $D=E\times_BE'$ be the fiber product of these bundles and let $\d_0D=\d_0E\times_BE'\cup E\times_B\d_0E'$. Then
\[
	\t(D,\d_0)=\chi(X,\d_0)\t(E,\d_0)+\chi(F,\d_0)\t(E',\d_0).\qed
\]
\end{cor}

\begin{xca}
Show that for any higher torsion invariant $\t$ and any unipotent
bundle pair $(E,\d_0)\to B$ we have
\[
    \t(E,\d_0)+(-1)^n\t(E,\d_1)=2\t^+(E,\d_0)
\]
where $n=\dim F$ is the fiber dimension.
\end{xca}


\subsection{Further extension}\label{subsection:further extension} 


Suppose that $(E,\d_0)\to B$ is a relatively unipotent smooth bundle pair. Then the torsion $\t(E,\d_0)$ can be defined as
follows.
\begin{enumerate}
    \item Let $D(\nu)\to E$ be the normal disk bundle and let
    $\d_0D=D(\nu)|\d_0E$.
    \item Embed $\d_0E$ fiberwise into the northern hemisphere of
    $B\times S^N$ for $N$ large.
    \item Thicken this embedding to a codimension zero embedding
    $\d_0D\to B\times S^N$. Then the union of $B\times D^{N+1}$ with $D$ along $\d_0D$ is
    unipotent.
    \item Define the higher torsion of $(D,\d_0)$ by
    \[
        \t(D,\d_0):=\t(B\times D^{N+1}\cup D).
    \]
    \item Define the higher torsion of $(E,\d_0)$ by
    \[
\t(E,\d_0):=\t(D,\d_0)-(s_1+s_2)tr^E_B(ch_{4k}(\nu))
    \]
    where $s_1,s_2$ are given in the explanation of the Main Theorem \ref{thm:tau even and tau odd are unique up to scalar}.
\end{enumerate}

The main theorem (and Theorem \ref{properties
of FR torsion}) are needed to show that this is well defined. The proof is that these formulas are well defined for $M_{2k}$ and for $\t_{2k}$ which span all possibilities by the main theorem. Therefore, we cannot use this extension to the relatively unipotent case to prove the main theorem.
In the case when
$(E,\d_0)$ is unipotent these formulas hold by additivity and
transfer.


\section{Stability of higher torsion} 


Smooth bundle are stabilized by taking products with disks. The following special case of the product formula says that higher torsion is a stable invariant.

\begin{cor}[stability of torsion]\label{cor:stability of torsion}
If $(E,\d_0)\to B$ is a unipotent smooth bundle pair then so is $(E\times D^n,\d_0E\times D^n)$ and the relative torsion is the same:
\[
	\t(E\times D^n,\d_0E\times D^n)=\t(E,\d_0).\qed
\]
\end{cor}

If $M$ is a compact smooth manifold then we recall that a \emph{concordance} of $M$ is a diffeomorphism of $M\times I$ which is the identity on $M\times0\cup \d M\times I$. Let $\cC(M)$ be the space of concordances of $M$ with the $C^\infty$ topology.
\[
	\cC(M)=\cD if\!f(M\times I\rel M\times0\cup \d M\times I)
\]
The classifying space $\cH(M)=B\cC(M)$ is the space of $h$-cobordisms of $M\rel \d M$. Recall that an \emph{$h$-cobordism} of $M\rel\d M$ is a compact smooth manifold $W$ with boundary $\d W=M\times 0\cup \d M\times I\cup M'$ where $M'$ is another compact smooth manifold with $\d M'=\d M\times I$. A mapping $B\to \cH(M)$ is a smooth bundle over $B$ whose fibers are all $h$-cobordisms of $M\rel \d M$ so that the subbundle with fiber $M\times 0\cup\d M\times I$ is trivial. We will call such a bundle an $h$-cobordism bundle over $B$.

There is a suspension map $\s:\cC(M)\to \cC(M\times I)$ which is highly connected when $\dim M$ is large by the concordance stability theorem (\cite{[I:Stability]} or the last chapter of \cite{[I:ComplexTorsion]}). The limit is the \emph{stable concordance space}
\[
	\cP(M)=\lim_\to \cC(M\times I^n)
\]
which is well known to be an infinite loop space. Therefore, the set $[B,B\cP(M)]$ of homotopy classes of maps from $B$ to the classifying space $B\cP(M)$ is an additive group. By the concordance stability theorem this group is isomorphic to $[B,\cH(M\times I^n)]$ for sufficiently large $n$.

\begin{prop}
For any $h$-cobordism bundle $E\to B$ let $\d_0E$ be the trivial subbundle with fiber $M\times 0\cup\d M\times I$. Then the higher torsion invariant $E\mapsto\t(E,\d_0)$ gives an additive map
\[
	\t:[B,\cH(M\times I^n)]\to\RR
\]
\end{prop}

\begin{proof} This is an immediate consequence of the horizontal additivity of $\t$ (Proposition \ref{horizontal additivity}) since the H-space structure on the $h$-cobordism space $\cH(M\times I^n)$ is given by lateral union. I.e., the sum of two mappings $B\to \cC(M\times I^n)$ is given by lateral union of the corresponding $h$-cobordism bundles followed by rescaling. (The lateral union is an $h$-cobordism of $M\times I^{n-1}\times [0,2]$. The last coordinate needs to be rescaled down to $[0,1]$.)
\end{proof}

\begin{cor}\label{trivial is same as rationally trivial}
Suppose that $E_\a$ is a collection of $h$-cobordism bundles which spans the $\QQ$ vector space $[B,\cH(M\times I^n)]\otimes\QQ$. Suppose also that $\t(E_\a,\d_0)=0$ for all $\a$. Then $\t(E,\d_0)=0$ for any $h$-cobordism bundle $E$ classified by a map $B\to \cH(M\times I^n)$.\qed
\end{cor}

\begin{remark}
This is equivalent to saying that, if $\t,\t'$ are two higher torsion invariants which agree on the rational generators $E_\a$ then they agree on all $h$-cobordism bundles $E$.
\end{remark}


\section{Computation of higher torsion} 


We will now show how the higher torsion invariants of unipotent bundles can be computed in many cases given the values of the parameters $s_1,s_2$. We recall that these parameters are given by
\[
    \t(S^n(\ll))=2s_nch_{4k}(\ll)
\]
where $S^n(\ll)$ is an $S^1$ or $S^2$ bundle associated to a complex line bundle $\ll$ over $B$.


\subsection{Torsion of disk and sphere bundles} 


\begin{theorem}
The higher torsion of the $D^n$-bundle $D^n(\xi)$ associated to an $SO(n)$-bundle $\xi$ over $B$ is given by
\[
    \t(D^n(\xi))=(s_1+s_2)ch_{4k}(\xi).
\]
\end{theorem}

\begin{proof}
For $n=2$ this is by definition of $\t$ in the boundary case:
\[
    \t(D^2(\xi))=\tfrac12\t(S^2(\xi))+\tfrac12\t(S^1(\xi))=(s_2+s_1)ch_{4k}(\xi).
\]
The general case follows from the product formula (Corollary \ref{cor:product formula}) and the splitting principle. If $n=2m+1$ then we may assume by the splitting principle that $\xi$ is a direct sum of $m$ complex line bundles $\ll_i$ and a trivial real line bundle. Then
\[
     \t(D^n(\xi))=\t\left(D^2(\ll_1)\times_B \cdots \times_B D^2(\ll_m)\times_B (B\times I)\right)
\]
\[
	=(s_1+s_2)\sum ch_{4k}(\ll_i)=(s_1+s_2)ch_{4k}(\xi)
\]
by the product formula. The even case is similar.
\end{proof}

From the calculation
\[
    \t(D^n(\xi))=\tfrac12\t(S^n)+\tfrac12\t(S^{n-1}(\xi))=(s_1+s_2)ch_{4k}(\xi)
\]
we get the following by induction on $n$.

\begin{cor} For $n>0$ the higher torsion of the $S^n$-bundle $S^n(\xi)$ associated an $SO(n+1)$-bundle $\xi$ over $B$ is given by
\[
    \t(S^n(\xi))=2s_nch_{4k}(\xi)
\]
where $s_n$ depend only on the parity of $n$. 
\end{cor}

Comparing this with Proposition \ref{MMM on sphere bundles} and Theorem \ref{FR torsion of linear sphere bundles} below we get the following.

\begin{lemma}\label{lem:tau is given on linear bundles}
Let $k\geq1$ and let $a,b\in\RR$ be given by
\[
    a:=\frac{(-1)^{k+1}2s_1}{\z(2k+1)},\quad b:=\frac{s_1+s_2}{(2k)!}
\]
Then $\t(E)=a\t_{2k}(E)+bM_{2k}(E)$ for $E$ any oriented linear sphere or disk bundle over $B$.
\end{lemma}

\begin{remark}
If the value of $\t$ on a bundle $E$ is determined by the known values of $\t$ on disk and sphere bundles then this lemma implies that $\t(E)=a\t_{2k}(E)+bM_{2k}(E)$ for that bundle.
\end{remark}


\subsection{Morse bundles} 


Given a Morse function $f:(M,\d_0M)\to (I,0)$, a compact $n$-manifold $M$ will be decomposed as a union of \emph{handles} $D^i\times D^{n-i}$ attached along $S^{i-1}\times D^{n-i}$ to the union of lower handles and the base $\d_0M\times [0,\e]$. Each such handle has a critical point at its center with index $i$. The \emph{core} of the handle is $D^i\times 0$. This is also the union of trajectories of the gradient of $f$ (with respect to some metric on $M$) which converge to the critical point. The tangent plane is the \emph{negative eigenspace} of the second derivative $D^2f$ at the critical point.

Suppose that $(E,\d_0)\to B$ is a smooth bundle pair and $f:E\to I=[0,1]$ is a fiberwise Morse function with $f^{-1}(0)=\d_0E$ and with distinct critical values. In other words, for each $b\in B$, the restriction $f_b:(F_b,\d_0)\to (I,0)$ is a Morse function with critical points $x_1(b),\cdots,x_m(b)$ having critical values $f_b(x_1)<f_b(x_2)<\cdots<f_b(x_m)$. 

For $1\leq j\leq m$ let $\xi_j$ be the negative eigenspace bundle of the fiberwise second derivative of $f$ along $x_j$. Let $\eta_j$ be the positive eigenspace bundle. Then $E$ has a filtration $E=E_m\supset E_{m-1}\supset\cdots\supset E_0$ where $E_0\cong \d_0E\times I$ and  
\[
	E_j=E_{j-1}\cup D(\xi_j)\times_B D(\eta_j)
\]
where $E_{j-1}\cap D(\xi_j)\times_B D(\eta_j)=S(\xi_j)\times_B D(\eta_j)$. By additivity we get the following.

\begin{lemma}
If $\d_0E\to B$ is unipotent and the bundles $\xi_j,\eta_j$ are oriented then $(E,\d_0)\to B$ is a unipotent bundle pair with torsion invariant
\[
	\t(E,\d_0)=\sum\t((D(\xi_j),S(\xi_j))\times_B D(\eta_j))
\]
\end{lemma}

Each summand in the above lemma can be determined using the product formula (Corollary \ref{cor:product formula}).

\begin{lemma}\label{lem:torsion of handle bundle} The value of $\t$ on the fiber product of an oriented linear disk bundle $D(\eta)$ and the oriented relative $i$-disk bundle $(D^i(\xi),S^{i-1}(\xi))$ is given by
\[
    \t((D^i(\xi),S^{i-1}(\xi))\times_BD(\eta))=(-1)^i\t(D(\eta))+\t(D^i(\xi),S^{i-1}(\xi))\]
    \[
    =(-1)^i(s_1+s_2)ch_{4k}(\eta)+(s_i-s_{i-1})ch_{4k}(\xi)
\]
\end{lemma}

\begin{remark}
Since $s_i-s_{i-1}=(-1)^is_2-(-1)^is_1$, the even and odd parts of the above formula are:
\[
    \t^+((D^i(\xi),S^{i-1}(\xi))\times_BD(\eta))=
    (-1)^is_2(ch_{4k}(\eta)+ch_{4k}(\xi))
\]
\[
    \t^-((D^i(\xi),S^{i-1}(\xi))\times_BD(\eta))=
    (-1)^is_1(ch_{4k}(\eta)-ch_{4k}(\xi))
\]
\end{remark}

Putting these together we get the following theorem which is a mild improvement over the obvious. Namely, the negative eigenspace bundles need not be oriented. However, the sum $\xi_i\oplus \eta_i$ of negative and positive eigenspace bundles must be orientable since it is the vertical tangent bundle along the $i$th component of the Morse critical set.

\begin{theorem}\label{torsion for Morse bundles}
Suppose that $(E,\d_0)\to B$ is a unipotent smooth bundle pair and $f:(E,\d_0E)\to (I,0)$ is a fiberwise Morse function with distinct critical values. Let $\xi_i, \eta_i$ be the negative and positive eigenspace bundles associated to the $i$th critical point. then
 \[
    \t(E,\d_0)=\sum(-1)^i(s_1+s_2)ch_{4k}(\eta)+(-1)^i(s_2-s_{1})ch_{4k}(\xi)
\]
\end{theorem}

\begin{proof}
Suppose first that the bundles $\xi_i,\eta_i$ are oriented. Then the two lemmas above apply to prove the theorem. If these bundles are not oriented then there is a finite covering $\wt{B}$ of $B$ so that, on the pull-back $\wt{E}$, the function $\wt{E} \to E\to I$ is a Morse function with oriented eigenspace bundles $\wt{xi}_i,\wt{\eta}_i$ which are pull-backs of $\xi_i,\eta_i$. Thus the theorem applies to $\wt{E}$. However, the induced map in real cohomology $H^\ast(B;\RR)\to H^\ast(\wt{B};\RR)$ is a monomorphism. The two expressions in our theorem are elements of $H^{4k}(B;\RR)$ which go to the same element of $H^{4k}(\wt{B};\RR)$. So, they must be equal.
\end{proof}


\subsection{Hatcher's example} 


One crucial example of a Morse bundle to which the above theorem holds is Hatcher's construction.  This constructs an exotic disk bundle over $B=S^n$ out of an element of the kernel of the $J$-homorphism
\[
    J:\pi_{n-1}O\to \pi_{n-1}^s(S^0).
\]
 B\"{o}kstedt \cite{[Bokstedt84]} interpretted this as a mapping from $G/O$ to the stable concordance space of a point $\cH(\ast)=B\cP(\ast)$. This means that a mapping $B\to G/O$ gives an exotic disk $h$-cobordism bundle over $B$.

Since $G/O$ is the homotopy fiber of the map $BO\to BG$, a map $B\to G/O$ is given by an $n$-plane bundle $\xi:B\to SO(n)$ together with a homotopy trivialization of the associated sphere bundle. I.e., we have a fiber homotopy equivalence
\[
    g:S^{n-1}(\xi)\simeq B\times S^{n-1}.
\]
We write this as a family of homotopy equivalences $g_t:S^{n-1}_t(\xi)\to S^{n-1}, t\in B$ and extend to the disk $D_t^n(\xi)$ (the fiber over $t\in B$ of the disk bundle $D^n(\xi)$)
\[
    \ol{g}_t:(D^n_t(\xi),S^{n-1}_t(\xi))\to (D^n,S^{n-1}).
\]
Assuming that $n$ and $m$ are large enough, we can lift $\ol{g}_t$ up to a family of embeddings
\[
    \wt{g}_t:(D^n_t(\xi),S^{n-1}_t(\xi))\to (D^n\times D^m,S^{n-1}\times D^m)
\]
If we let $\eta$ be the complementary bundle to $\xi$, we can extend this to a family of codimension zero embeddings
\[
    G_t:(D^n_t(\xi)\times D^m_t(\eta),S^{n-1}_t(\xi)\times D^m_t(\eta))
    \to
     (D^n\times D^m,S^{n-1}\times D^m)
\]

We will use this to construct a nontrivial $h$-cobordism bundle over $B$. We start this construction with the trivial $h$-cobordism bundle $E_0=B\times D^{n+m-1}\times I$. Over each $t\in B$ the fiber is $D^{n+m-1}\times I$ which the trivial $h$-cobordism of $D^{n+m-1}$. We now add two canceling handles: an $n-1$ handle $H^{n-1}$ and an $n$ handle $H^n$. We will then replace $H^n$ with a variable $n$-handle $H_t^n$.

The $n-1$ handle $H^{n-1}$ is simply $S^{n-1}\times D^{m+1}$ which is attached by boundary connected sum to the top boundary $D^{n+m-1}\times 1$. The next cell $H^n$ is a copy of $D^n\times D^m$ which is attached along an embedding:
\[
	id\times \a:S^{n-1}\times D^m\to S^{n-1}\times \d D^{m-1}
\]
where $\a:D^m\to \d D^{m-1}$ is some fixed embedding.

The new $h$-cobordism bundle $B\times(D^{n+m-1}\times I\cup H^{n-1}\cup H^n)$ is still a trivial bundle. However, it contains a nontrivial subbundle given by replacing $H^n$ with the image of $G_t$ for every $t\in B$. Let $H_t^n$ denote this subset. Then $H_t^n$ is a variable (``rotating'') $n$-handle which is attached to $\d H^{n-1}$ by the mapping
\[
	(id\times\a)\circ G_t:S^{n-1}_t(\xi)\times D^m_t(\eta)\to S^{n-1}\times \d D^{m-1}\subset\d H^{n-1}
\]
Let $(E_1,\d_0)$ denote the bundle pair with fiber pair
\[
	(D^{n+m-1}\times I\cup H^{n-1}\cup H_t^n, D^{n+m-1}\times 0\cup \d D^{n+m-1}\times I).
\]
The main property of this bundle are given by the following theorem originally due to B\"{o}ckstedt \cite{[Bokstedt84]} (see also \cite{[Bokstedt-Waldhausen87]}) and reproved using Morse theory in \cite{[I:BookOne]}.

\begin{theorem}\label{Bockstedt}
The bundles $(E_1,\d_0)\to B$ given by Hatcher's construction above are rational generators for the group $[B,\cH(D^{n+m-1})]$
\end{theorem}

By construction, Hatcher's $h$-cobordism bundle admits a fiberwise Morse function with exactly two critical points on each fiber. The first critical point has index $n-1$ and has a trivial negative eigenspace bundle. The second has index $n$ and its negative eigenspace bundle is $\xi$. By Theorem \ref{torsion for Morse bundles}, the torsion of Hatcher's bundle is:
\[
    \t(\Delta)(\xi)=(-1)^n(s_1+s_2)ch_{4k}(\eta)+(-1)^n(s_2-s_1)ch_{4k}(\xi).
\]
Since $ch_{4k}(\eta)+ch_{4k}(\xi)=0$ ($\xi\oplus\eta$ being trivial) this simplifies to:
\begin{equation}\label{torsion of Hatcher}
    \t(\Delta)(\xi)=(-1)^{n+1}2s_1ch_{4k}(\xi).
\end{equation}
Since $s_1=0$ for $M_{2k}$ this implies the following.

\begin{lemma}
For any bundle $(D,\d_0)\to B$ of $h$-cobordisms of disks we have $\t(D,\d_0)=\t(D)=a\t_{2k}(D)$ where $a$ is given in Lemma \ref{lem:tau is given on linear bundles}.
\end{lemma}

\begin{proof}
This follows from the calculation (\ref{torsion of Hatcher}) above in the case of Hatcher's example. By B\"{o}ckstedt (Theorem \ref{Bockstedt}) and Corollary \ref{trivial is same as rationally trivial} this implies that $\t(D)=a\t_{2k}(D)$ for any disk $h$-cobordism bundle of sufficiently large fiber dimension.

If $D$ is any disk $h$-cobordism bundle then $\t(D)=\t(D\times D^k)$ and $\t_{2k}(D)=\t_{2k}(D\times D^k)$ by stability (Corollary \ref{cor:stability of torsion}). Furthermore, $D\times D^k$ will be a disk bundle with fiber dimension, say $N$, which contains a trivial $N-1$ disk boundary in its vertical boundary. Therefore it can be modified to be an $h$-cobordism bundle where the formula holds. It follows from the stable result that $\t(D\times D^k)=a\t_{2k}(D\times D^k)$. The lemma follows.
\end{proof}

\begin{prop}\label{all disk bundles are good}
For any smooth oriented disk bundle $D\to B$ we have $\t(D)=a\t_{2k}(D)+bM_{2k}(D)$ where $a,b$ are given in Lemma \ref{lem:tau is given on linear bundles}.
\end{prop}

\begin{proof}
By Corollary \ref{cor:stability of torsion}, we can stabilize $D$ without changing the value of $\t(D),\t_{2k}(D)$ or $M_{2k}(D)$. Then the fiber of $D\to B$ will be a disk $D^N$ where $N>\dim B$. Next, choose a smooth section $B\to \d^vD$. Thicken this up to get a linear $N-1$ disk bundle $D(\xi)$ which lies in $\d^vD$. 

Suppose first that $\xi$ is a trivial bundle. Then the bundle $D$ is equivalent to an $h$-cobordism bundle. We can see this by attaching $D$ along the image of $D(\xi)$ to the top of $B\times D^{N-1}\times I$. Then, by additivity, the resulting $h$-cobordism bundle will have the same torsion ($\t,\t_{2k}$ and $M_{2k}$) as the disk bundle. The equation relating these torsions holds for the $h$-cobordism bundle by the lemma. Therefore it holds for the disk bundle.

If the linear bundle $\xi$ is nontrivial we simply take the fiber product of $D$ with the disk bundle of the complementary bundle $\eta$. The bundles $\xi,\eta$ are necessarily oriented since $D$ is oriented. So, the new bundle $D\times_BD(\eta)$ has torsion
\[
    \t(D\times_B D(\eta))=\t(D)+\t(D(\eta))
\]
by the product formula. Since the proposition holds for $D\times_B D(\eta)$ and $D(\eta)$, it holds for $D$.
\end{proof}


\section{Uniqueness of higher torsion} 


We are now ready to prove the uniqueness part of the main theorem. Namely, we will show that $\t=a\t_{2k}+bM_{2k}$ on all unipotent bundles. To do this we show that the difference is zero. We continue to assume that $\t_{2k},M_{2k}$ are higher torsion invariants, facts that we prove in the last section.


\subsection{The difference torsion} 


We define the \emph{difference torsion} by
\[
    \t^\delta:=\t-a\t_{2k}-bM_{2k}.
\]
This is a linear combination of higher torsion invariants and
therefore a higher torsion invariant. Furthermore this new invariant has
the property that $s_1=s_2=0$. So:

\begin{lemma}\label{tdelta=0 on disk and sphere bundles}
The difference torsion is zero on all smooth oriented disk bundles and oriented linear
sphere bundles.
\end{lemma}

\begin{lemma}[thickening lemma]\label{tdelta is the same on linear disk bundles}
If $q:D\to E$ is an oriented linear disk bundle then the
difference torsion of any unipotent bundle pair $(E,\d_0)$ is
equal to that of $(D,\d_0)$ as a bundle over $B$ where
$\d_0D=q^{-1}(\d_0E))$. I.e.,
\[
    \t_B^\delta(D,\d_0)=\t^\delta(E,\d_0).
\]
\end{lemma}

\begin{proof}
By the relative transfer formula we have
\[
    \t^\delta_B(D,\d_0)
    =\t^\delta(E,\d_0)+tr^E_B(\t^\delta_E(D)).
\]
But $\t^\delta_E(D)=0$ since $D$ is a linear disk bundle over $E$.
\end{proof}

\begin{lemma}\label{tdelta=0 on H-cobordisms}
The difference torsion is zero on any unipotent bundle pair
$(E,\d_0)\to B$ whose fibers $(F,\d_0)$ are $1$-connected and homologically trivial ($H_\ast(F,\d_0F;\ZZ)=0$).
\end{lemma}

\begin{proof} First we note that the homological condition on the fiber implies that $F/\d_0F$ is contractible (being $1$-connected and acyclic).

Choose a smooth fibred embedding 
\[
	(E,\d_0E)\to(B\times S^{N-1}\times I,B\times S^{N-1}\times 0)
\]
for some large $N$. Let $(D,\d_0)$ be a tubular neighborhood of the image of $(E,\d_0)$. Since $(D,\d_0)$ is the normal disk bundle it has the same difference torsion as $(E,\d_0)$. So it suffices to show that $\t^\delta(D,\d_0)=0$.

The union $D\cup B\times D^N$ is a bundle whose fiber is homotopy equivalent to $F/\d_0F$ and thus contractible. Therefore it is a disk bundle when the corners are rounded. So, its difference torsion is trivial. The additivity theorem implies that
\[
	\t^\delta(D\cup B\times D^N)=\t^\delta(D)+\t^\delta(B\times D^N)-\t^\delta(\d_0D)=\t^\delta(D,\d_0)=0
\]
as claimed.
\end{proof}

\begin{lemma}[main lemma]\label{main lemma}
The difference torsion $\t^\delta$ is a fiber homotopy invariant
of unipotent smooth bundle pairs.
\end{lemma}

\begin{proof}
Suppose that $(E_1,\d_0)$ and $(E_2,\d_0)$ are unipotent smooth
bundle pairs over $B$ that are fiber homotopy equivalent. Then we
want to show that
\[
    \t^\delta(E_1,\d_0)=\t^\delta(E_2,\d_0).
\]
By the thickening lemma we can make the fiber dimension of $(E_2,\d_0)$ arbitrarily large. Then we can approximate the fiber homotopy equivalence by
a fiberwise smooth embedding
\[
    g:(E_1,\d_0)\to(D_2,\d_0).
\]
Using the thickening lemma again we can assume that $g$ is a codimension $0$ embedding.
The complement of the image of $g$ is a unipotent fiberwise
$h$-cobordism with trivial $\t^\delta$ by the previous lemma. Thus
$\t^\delta(E_1,\d_0)=\t^\delta(D_2,\d_0)$ by additivity.
\end{proof}


\subsection{Vanishing of the fiber homotopy invariant} 


The final step in the proof of Theorem \ref{thm:tau even and tau odd are unique up to scalar} is to show the following.

\begin{lemma}\label{lem:tdelta=0 relative}
Any higher torsion invariant $\t^\delta$ which is also a fiber
homotopy invariant of unipotent smooth bundle pairs must be zero.
\end{lemma}

Since $\t^\delta$ is a fiber homotopy invariant, it is well
defined on any fibration pair $(Z,C)\to B $ with fiber $(X,A)$
which is \emph{unipotent} in the sense that $H_\ast(X;\QQ)$ and
$H_\ast(A;\QQ)$ are unipotent as $\pi_1B$-modules and
\emph{smoothable} in the sense that it is fiber homotopy
equivalent to a smooth bundle pair $(E,\d_0)$ with compact
manifold fiber $(F,\d_0)$. By definition, the relative torsion is related to the absolute
torsion by
\begin{equation}\label{eq:def of relative difference torsion}
    \t^\delta(Z,C)=\t^\delta(Z)-\t^\delta(C).
\end{equation}

We will examine what it means for a fibration $X\to Z\to B$ to be unipotent and smoothable in terms of diffeomorphism spaces of compact manifolds. Suppose that $M\to E\to B$ is a smoothing of $Z\to B$. Then $M$ is a compact smooth manifold homotopy equivalent to $X$. By taking the product with a high dimensional disk we may assume that the inclusion map of the boundary $\d M\hookrightarrow M$ is highly (at least $4k$) connected. The structure group of the bundle $E\to B$ will be a subgroup of the diffeomorphism group of $M$ consisting of orientation preserving diffeomorphisms which act unipotently with respect to a fixed flag $\cF$ in the rational homology of $M$, i.e., these are diffeomorphisms which preserve the flag and act as the identity on the successive quotients. Let $Dif\!f_\cF(M)$ denote the space of such diffeomorphisms with the $C^\infty$ topology. Then the smooth bundle $E\to B$ is classified by a mapping $B\to BDif\!f_\cF(M)$.

The difference torsion comes from a universal invariant
\[
	\t^\delta\in H^{4k}(BDif\!f_\cF(M);\RR)
\]
associated to the universal bundle over $BDif\!f_\cF(M)$ with fiber $M$. Let $M_0\subset M$ be the closure of the complement of a tubular neighborhood of $\d M$. Then $M_0$ is a deformation retract of $M$. So, any diffeomorphism of $M$ which is the identity on $M_0$ will automatically induce the identity in homology. Thus we have an inclusion map $j:Dif\!f(M\rel M_0)\hookrightarrow Dif\!f_\cF(M)$.

\begin{prop}\label{pull back of universal tdelta is zero}
The universal difference invariant $\t^\delta$ maps to zero under the induced map
\[
 	j^\ast:H^{4k}(BDif\!f_\cF(M);\RR)\to H^{4k}(BDif\!f(M\rel M_0);\RR).
\]
\end{prop}

\begin{proof}
When the universal bundle with fiber $M$ is pulled back to $BDif\!f(M\rel M_0)$ it will contain a trivial $M_0$ bundle. The remainder is an $h$-cobordism bundle for $\d M$. Therefore $\t^\delta$ of this pull-back is zero by additivity and Lemma \ref{tdelta=0 on H-cobordisms}.
\end{proof}

\begin{prop}\label{prop: fibration sequence for diff M}
There is a fibration sequence
\[
	Dif\!f(M\rel M_0)\to Dif\!f_\cF(M)\xrarrow{r} Emb_\cF(M_0,int M)
\]
where $Emb_\cF(M_0,int M)$ is the space of orientation preserving smooth embeddings $M_0\to int M$ which is unipotent on $H_\ast(M_0;\QQ)\cong H_\ast(M;\QQ)$ with respect to the flag $\cF$ and $r$ is the restriction map.
\end{prop}

\begin{proof}
See the appendix of \cite{[Cerf:gamma-four]}.
\end{proof}

\begin{lemma}\label{lem:MV for tdelta} Given unipotent smoothable
fibrations $(Z,C)$ and $Y$ over $B$ and any continuous mapping
$f:C\to Y$ over $B$ the union $Y\cup_CZ$ is also unipotent and
smoothable with
\[
    \t^\delta(Y\cup_CZ)=\t^\delta(Y)+\t^\delta(Z,C).
\]
\end{lemma}

\begin{proof}
Unipotence follows from the Mayer-Vietoris sequence for the fiber
homology of $Y\cup_CZ$. Additivity of torsion gives the formula.
Smoothing is easy, very similar to the construction in subsection \ref{subsection:further extension}.
\end{proof}

\begin{lemma}\label{lem:fiberwise suspension is smoothable}
For any unipotent smoothable fibration pair $(Z,C)$ the fiberwise
smash $Z/\!_BC$ and fiberwise suspension $\Sig_B(Z/\!_BC)\to B$
(with fibers $X/A$ and $\Sig(X/A)$) are unipotent and smoothable
with $
    \t^\delta(Z/\!_BC)=\t^\delta(Z,C)=-\t^\delta(\Sig_B(Z/\!_BC)).
$
\end{lemma}

\begin{proof} Taking $Y=B$ in the previous lemma we get 
\[
	\t^\delta(Z/\!_BC)=\t^\delta(B\cup_CZ)=\t^\delta(Z,C)=\t^\delta(Z)-\t^\delta(C).
\]
Apply this to the cofibration sequence $X/A\vee X/A\to X/A\to
\Sig (X/A)$ to get
\[
    \t^\delta(\Sig_B(Z/\!_BC))=\t^\delta(Z/\!_BC)-2\t^\delta(Z/\!_BC)
    =-\t^\delta(Z/\!_BC).
\]
Where we use the equation $\t^\delta(X\vee Y)=\t^\delta(X)+\t^\delta(Y)$ which is another special case of Lemma \ref{lem:MV for tdelta}.
\end{proof}

The following lemma and its proof goes along the lines of a
discussion I had with John R. Klein.

\begin{lemma}\label{lem:tdelta=0 for rationally trivial fibers}
$\t^\delta(Z)=0$ for pointed unipotent smoothable fibrations with
rationally acyclic fibers.
\end{lemma}

\begin{proof} By \emph{pointed} we mean that $Z$ has a section
$B\to Z$ whose image is a copy of $B$. By the previous lemma we
may assume that the fiber $X$ has been suspended many more times
than the dimension of $B$. Since $X$ is a rationally trivial finite complex it has only finitely many homotopy classes of self maps. Thus there is a finite covering $\wt{B}$ of $B$ so that the action of $\pi_1\wt{B}$ on $X$ is trivial up to pointed homotopy. Since $H^{4k}(B,\RR)$ maps monomorphically into $H^{4k}(\wt{B};\RR)$ it suffices to show that the difference torsion of the pull back $\wt{Z}$ of $Z$ to $\wt{B}$ is trivial.

By assumption there is a smooth thickening $M\to E\to \wt{B}$ of the bundle $\wt{Z}\to \wt{B}$ where $M$ is a compact $n$-manifold homotopy equivalent to $X$. By passing to the normal disk bundle we may assume that the vertical tangent bundle of $E$ is trivial. Since the bundle has a section this means that $E$ contains a trivial $n$-disk bundle $\wt{B}\times D^n$. Therefore the bundle is classified by a mapping $\wt{B}\to BDif\!f_0(M\rel D^n)$ where $Dif\!f_0$ means diffeomorphisms tangentially homotopic to the identity.\vs2

\ul{Claim 1}: We may assume that $Dif\!f_0(M\rel D^n)$ is connected and therefore $BDif\!f_0(M\rel D^n)$ is simply connected.\vs2

Proof: Let $M_0$ be the closure of the complement of a collar neighborhood of $\d M$ in $M$. We may assume the disk $D^n$ lies in $M_0$. By taking the product with a disk if necessary we may assume that $M_0$ contains as a deformation retract a spine $K\simeq X$ of dimension $<(n-4k)/2$. Let $f$ be any diffeomorphism of $M$ which is tangentially homotopic to the identity. By immersion theory, the restriction of $f$ to $M_0$ is isotopic by regular immersions to the inclusion map. Since $M_0$ has a small dimensional spine we can deform this immersion isotopy into an isotopy through embeddings. This extends to $M$ by isotopy extension. But $\d M$ is simply connected. So, its concordance space $\cC(\d M)$ is connected by Cerf \cite{[Cerf]}. Thus $f$ is isotopic to the identity as claimed.
\vs2

\ul{Claim 2}: There is a fibration sequence:
\[
	Emb_0(M_0,int M\rel D^n)\to BDif\!f(M\rel M_0)\xrarrow{Bj} BDif\!f_0(M\rel D^n) 
\]
where $Emb_0$ means embeddings isotopic to the inclusion map.\vs2

Proof: Choose a fixed embedding $\f_0:M\to \RR^\infty$ and let $EDif\!f(M\rel D^n)$ be the space of all embeddings $\f:M\to \RR^\infty$ which agree with $\f_0$ on $D^n$. Then $EDif\!f(M\rel D^n)$ is contractible and admits a free action by $Dif\!f_0(M\rel D^n)$. The quotient space $BDif\!f_0(M\rel D^n)$ is the space of all pairs $(W,[\f])$ where $W$ is a submanifold of $\RR^\infty$ containing $D^n$ and $[\f]$ is an isotopy class of diffeomorphisms $\f:M\to W$ which are equal to $\f_0$ on $D^n$. The classifying space for $Dif\!f(M\rel M_0)$,
\[
	BDif\!f(M\rel M_0)=EDif\!f(M\rel D^n)/Dif\!f(M\rel M_0),
\]
is the space of all pairs $(W,\psi)$ where $W\subset\RR^\infty$ is as above and $\psi:M_0\to W$ is an embedding equal to $\f_0$ on $D^n$ and whose image has complement a collar neighborhood of $\d W$. Each $\psi$ determines an isotopy class $[\f]$ and, for $W=M$, the space of all $\psi$ which give the isotopy class $[id]$ of the identity is exactly $Emb_0(M_0,int M\rel D^n)$ as claimed.
\vs2

\ul{Claim 3}: $Emb_0(M_0,int M\rel D^n)$ is rationally trivial through degree $4k$.\vs2

Proof: By immersion theory and transversality this space has the $4k$ homotopy type of the identity component of the space of all pointed maps $X\to X\times O$ which is rationally trivial since $X$ is rationally trivial.

One consequence of this is that the mapping $Dj$ in Claim 2 is a rational homotopy equivalence through degree $4k$. Thus $Dj$ induces an isomorphism in rational cohomology in degee $4k$. But $Dj^\ast(\t^\delta)=0$ by Proposition \ref{pull back of universal tdelta is zero}. Therefore, the universal class $\t^\delta$ is trivial.
\end{proof}

The rest of the proof of the Lemma \ref{lem:tdelta=0 relative}
follows suggestions of E. Bruce Williams.

\begin{lemma}\label{lem:tdelta=0}
$\t^\delta(Z)=0$ on all unipotent smoothable fibration.
\end{lemma}

\begin{proof} Taking the fiberwise suspension we may assume that
the fibers are pointed. By Lemma \ref{lem:tdelta=0 for rationally
trivial fibers} it suffices to reduce the rank of the rational
homology of the fiber $X$ of a pointed smoothable fibration
without changing the value of $\t^\delta(Z)$. We use the fact that stable homotopy groups are rationally the same as reduced homology.

Let $n$ be maximal
so that the reduced rational homology $\ol{H}_n(X;\QQ)$ is
nonzero. Choose one generator which is fixed by the action of
$\pi_1B$. (This exists since the action is unipotent.) Assuming that $X$ has been suspended a large number of
times, a multiple of this generator is represented by some
$\a\in\pi_n(X)$. Multiplying by the order of the torsion subgroup
of $\pi_n(X)$ we may assume that $\a$ represents a rationally
nontrivial element of $\pi_n(X)$ which is fixed by the action of
$\pi_1B$.

The homotopy groups $\pi_i(X)$ will be torsion and thus finite for
all $n<i< n+\dim B$. Let $m$ be the product of the orders of all
of these groups. Then, by obstruction theory, $m\a$ is represented
by a fiber preserving continuous mapping
\[
    f:(B\times S^n,B\times\ast)\to (Z,B).
\]
The fiberwise mapping cone $C_B(f)$ of $f$ has the same torsion:
$\t^\delta(C_B(f))=\t^\delta(Z)$ by Lemma \ref{lem:MV for tdelta}
and its fiber $X\cup e^{n+1}$ has less rational homology than $X$. The lemma follows by induction on the rank of the rational homology of $X$.
\end{proof}

This completes the proof of Lemma \ref{lem:tdelta=0 relative} which implies the Main Theorem \ref{thm:tau even and tau odd are unique up to scalar}.


\section{Existence of higher torsion}\label{sec:existence} 


In this section we show that higher Miller-Morita-Mumford classes $M_{2k}$ and higher Franz-Reidemeister torsion  $\t_{2k}$ are linearly independent higher torsion invariants.


\subsection{Miller-Morita-Mumford classes} 


If $p:(E,\d_0)\to B$ is any smooth bundle pair (not necessarily unipotent) with compact fiber $(F,\d_0)$, the (higher relative) \emph{Miller-Morita-Mumford} class $M_{2k}$ is defined to be the integral cohomology classes given by
\[
    M_{2k}(E,\d_0):=tr^{(E,\d_0)}_B((2k)!ch_{4k}(T^vE))
\in H^{4k}(B;\ZZ)
\]
for $k\geq 1$ where $T^vE$ is the vertical tangent bundle of $E$, $ch_{4k}(\xi):=\frac12 ch_{4k}(\xi\otimes\CC)$ is the \emph{Chern character} and $tr^{(E,\d_0)}_B$ is the relative transfer defined in (\ref{relative transfer}). In degree 0 this formula gives a half integer:
\[
    M_0(E,\d_0)=\frac{n}2\chi(F,\d_0)=\frac{n}2(\chi(F)-\chi(\d_0F))
\]
where $n=\dim F$. The properties of the transfer \cite{[Becker-Gottlieb:Adams-conj]} give the following properties of the higher Miller-Morita-Mumford classes.

\begin{lemma}[vertical additivity of $M_{2k}$]\label{vertical additivity of MMM}
Suppose that $E$ is the union of two subbundle $E=E_1\cup E_2$ and $E_1\cap E_2=\d_0E_2\subseteq \d_1E_1$. Then
\[
	M_{2k}(E,\d_0E_1)=M_{2k}(E_1,\d_0)+M_{2k}(E_2,\d_0).
\]
\end{lemma}

\begin{proof}
Choose a vertical tangent vector field for $E$ which points inward along $\d_0E_1$, outward along the rest of $\d^vE$ and points from $E_1$ to $E_2$ along $\d_0E_2$. Then the zero set of this vector field, say $Z$, is a disjoint union $Z= Z_1\coprod Z_2$ where $Z_i\subset E_i-\d^vE_i$. $M_{2k}(E,\d_0E_1)$ is given by restricting a multiple of the Chern character of $T^vE$ to $Z$ and pushing down to $B$. (See \cite{[I:ComplexTorsion]} for a detailed discussion of this.) This is obviously the sum of the push-downs of the restrictions to $Z_1$ and $Z_2$ which give $M_{2k}(E_1,\d_0)$ and $M_{2k}(E_2,\d_1)$.
\end{proof}

\begin{prop}\label{properties of MMM} For all $k>0$ we have
\begin{enumerate}
    \item \emph{(stability)} $M_{2k}(E\times I)=M_{2k}(E)$.
    \item \emph{(relative formula)} $M_{2k}(E,\d_0)=M_{2k}(E)-M_{2k}(\d_0E)$.
    \item \emph{(additivity)} If $E_1,E_2$ are smooth bundles over $B$ with the same vertical boundary $E_1\cap E_2=\d^vE_1=\d^vE_2$ then
    \[
        M_{2k}(E_1\cup E_2)=M_{2k}(E_1)+M_{2k}(E_2)-M_{2k}(\d^vE_1).
    \]
    \item \emph{(transfer)} If $q:D\to E$ is a bundle with fiber $X$ then $M_{2k}(D)_B\in H^{4k}(B)$ and $M_{2k}(D)_E\in H^{4k}(E)$ are related by
    \[
        M_{2k}(D)_B=\chi(X)M_{2k}(E)+tr^E_B(M_{2k}(D)_E).
    \]
\end{enumerate}
\end{prop}

\begin{remark}
The relative formula and additivity in the case $E_1\cong E_2$ are important since they imply that the relative invariant $M_{2k}(E,\d_0)$ is related to the closed fiber case in the way that higher torsion invariants are supposed to behave.
\end{remark}

\begin{proof} The transfer formula follows from the formula $T^v_BD=T^v_ED\oplus q^\ast T^vE$ and the additivity of the Chern character:
\begin{align*}
    M_{2k}(D)_B&=tr^D_B((2k)!ch_{4k}T^v_BD)\\
    &=tr^E_Btr^D_E((2k)!ch_{4k}(T^v_ED))
    +tr^E_Btr^D_E(q^\ast(2k)!ch_{4k}(T^vE))\\
    &=tr^E_B(M_{2k}(D)_E)+\chi(X)M_{2k}(E)
\end{align*}
Stability is a special case of the transfer formula. 

Vertical additivity (Lemma \ref{vertical additivity of MMM}) applied to $E\cong E\cup \d_0E\times I$ gives
\[
	M_{2k}(E,\d_0)=M_{2k}(E)-M_{2k}(\d_0E\times I)
\]
By stability this gives the relative formula.

Vertical additivity applied to $E=E_1\cup E_2$ gives
\[
	M_{2k}(E)=M_{2k}(E_1,\d^v)+M_{2k}(E_2)
\]
This gives additivity by the relative formula.
\end{proof}

\begin{lemma}\label{MMM on sphere bundles}
If $S^{2n}(\xi)\to B$ is the $S^{2n}$-bundle associated to an $SO(2n+1)$-bundle $\xi$ over $B$ then
\[
    M_{2k}(S^{2n}(\xi))=2(2k)!ch_{4k}(\xi).
\]
\end{lemma}

\begin{proof}
This follows from additivity and the elementary calculation
\[
    M_{2k}(D^m(\xi))=(2k)!ch_{4k}(\xi)
\]
for any linear disk bundle $D^m(\xi)$.
\end{proof}

\begin{theorem}\label{MMM is an even torsion invariant}
The higher Miller-Morita-Mumford class $M_{2k}$ is a nontrivial even higher torsion invariant for all $k>0$ with $s_1=0$ and $s_2=(2k)!$.
\end{theorem}

\begin{proof}
The additivity axiom follows from additivity (Proposition \ref{properties of MMM}.3) applied to $E_1\cup E_2$, $DE_1$ and $DE_2$. The transfer axiom is a special case of the transfer formula (Proposition \ref{properties of MMM}.4). The calculation of $s_2$ follows from the above lemma. Finally, $s_1=0$ since $e(E)=(-1)^ne(E)$ for an oriented bundle with closed $n$-dimensional fiber.
\end{proof}


\subsection{Higher FR-torsion} 


Higher \emph{Franz-Reidemeister} (FR)-torsion invariants are real characteristic class
\[
    \t_{2k}(E,\d_0)\in H^{4k}(B;\RR)
\]
for $k\geq1$ defined for smooth bundle pairs
$(F,\d_0)\to(E,\d_0)\to B$ which are \emph{relatively unipotent}
in the sense that $H_\ast(F,\d_0F;\QQ)$ is unipotent as a $\pi_1B$
module. (See \cite{[I:BookOne]}.) The main tool for computing the higher FR-torsion is the \emph{framing principle} which is proved in a special case in \cite{[I:BookOne]} and in general in \cite{[I:ComplexTorsion]}.

We will state the framing principle leaving unanswered the difficult question: What exactly is a ``family of chain complexes'' and how is the higher torsion of such a family defined? Suppose that $(E,\d_0)\to B$ is a relatively unipotent smooth bundle pair and $f:E\to\RR$ is a fiberwise oriented generalized Morse function (i.e., it has Morse and birth-death critical point and an orientation of the negative eigenspace of $D^2f_t$ along the Morse point sets so that they cancel with positive incidence along the birth-death points). Suppose that the vertical gradient $\nabla^vf$ of $f$ with respect to some metric points inward along $\d_0E$ and outward along $\d_1E$ and point from $\d_0E$ towards $\d_1E$ along $\d_0E\cap \d_1E$ (i.e., the zero set $\Sig$ of $\nabla^vf$ is dual to the relative Euler class $e(E,\d_0)$).

\begin{theorem}[Framing Principle \cite{[I:ComplexTorsion]}]
Let $C(f_t), t\in B$ be the family of cellular chain complexes obtained from the functions $f_t$. Then
\[
	\t_{2k}(E,\d_0)=\t_{2k}(C(f_t))+(-1)^k\z(2k+1)p^\Sig_\ast\left(
		ch_{2k}(\xi)
	\right)
\]
Here $\xi$ is the negative eigenspace bundle of $D^2f_t$ along the singular set, the push-down is the alternating sum $p_\ast^\Sig=\sum(-1)^ip_\ast^i$ where $p_\ast^i$ is the restriction of the push-down on the set $\Sig^i(f_t)$ of critical points of index $i$.
\end{theorem}

\begin{remark}
The higher FR-torsion of a family of based free chain complexes $C(f_t)$ is defined provided that its rational homology is unipotent. We call $\t_{2k}(C(f_t))$ the \emph{algebraic torsion}. The other summand is called the \emph{correction term}.
\end{remark}

An easy example of this rule is given by a linear $n$-sphere bundle $S^n(\xi)$ associated with an oriented $n$-plane bundle $\xi$ over $B$. This has a fiberwise Morse function $f$ having exactly two critical point in each fiber of index $0$ and $n$. The family of chain complexes $C(f_t)$ is constant and therefore has trivial algebraic torsion. The negative eigenspace bundle is trivial over $\Sig^0(f_t)$ and isomorphic to $\xi$ over $\Sig^n(f_t)$. Therefore the framing principle gives the following.

\begin{cor}\label{cor:FR torsion for sphere bundle with section} The higher FR-torsion of an oriented linear sphere bundle is given by
\[
	\t_{2k}(S^n(\xi))=(-1)^{n+k}\z(2k+1) ch_{2k}(\xi)
\]
if $\xi$ is an oriented $n$-plane bundle.\qed
\end{cor}

\begin{theorem}\cite{[I:BookOne]}\label{properties of FR torsion} The higher FR-torsion $\t_{2k}(E,\d_0)\in H^{4k}(B;\RR)$ is defined for relatively unipotent smooth bundle pairs $(E,\d_0)\to B$ and satisfies the following conditions.
\begin{enumerate}
\item \emph{(relative formula)} If $E$ and $\d_0E$ are unipotent bundles then 
\[
	\t_{2k}(E,\d_0)=\t_{2k}(E)-\t_{2k}(\d_0E).
\]
    \item \emph{(additivity)} If $(E_1,\d_0),(E_2,\d_0)$ are relatively unipotent smooth bundles over $B$ with $E_1\cap E_2=\d_0E_2\subseteq \d_1E_1$ then
    \[
        \t_{2k}(E_1\cup E_2,\d_0E_1)=\t_{2k}(E_1,\d_0)+\t_{2k}(E_2,\d_0).   \]
    \item \emph{(stability)} If $q:D\to E$ is a linear disk bundle and $\d_0D=q^{-1}(\d_0E)$ then the higher FR-torsion of $(D,\d_0)$ as a bundle pair over $B$ is equal to the higher torsion of $(E,\d_0)$:
    \[
        \t_{2k}(D,\d_0)_B=\t_{2k}(E,\d_0).
    \]
\end{enumerate}
\end{theorem}

\begin{proof}
The proofs of these theorems, given in great detail in \cite{[I:BookOne]}, can be summarized as follows. For additivity we may assume that there is a fiberwise oriented generalized Morse function (GMF) $f$ on $(E,\d_0)$ whose restriction to $(E_1,\d_0)$ and $(E_2,\d_0)$ is suitable for defining their torsion. Then the family of chain complexes $C(f_t)$ fits into a short exact sequence
\[
	0\to C(f|E_1)\to C(f_t)\to C(f|E_2)\to 0
\]
where $C(f|E_i)$ is the family of chain complexes associated to the restriction of $f$ to $E_i$. However, algebraic torsion is additive for short exact sequences, i.e., $\t_{2k}(C(f_t))=\t_{2k}(f|E_1)+\t_{2k}(f|E_2)$. The correction term is also additive since the critical set of $f_t$ will be a disjoint union of the critical sets of $f|E_1$ and $f|E_2$.

For stability we note that a fiberwise oriented GMF $f$ for $E$ gives a fiberwise oriented GMF $\s f$ on any linear disk bundle $p:D\to E$ by $\s f(x)=f(p(x))+\length{x}^2$. This has the same cellular chain complex and the same critical set with the same negative eigenspace bundles. So, the torsion is unchanged.

The relative formula follows formally from additivity in the case $E\cong E\cup \d_0E\times I$ and stability which gives $\t_{2k}(\d_0E\times I)=\t_{2k}(\d_0E)$.
\end{proof}

Higher FR-torsion and the higher Miller-Morita-Mumford classes are related by:

\begin{theorem}\cite{[I:ComplexTorsion]}\label{even part of FR is MMM}
If $E\to B$ is a unipotent smooth bundle with closed even
dimensional fibers and $k\geq1$ then
\[
    \t_{2k}(E)
    =(-1)^k\frac{\z(2k+1)}{2(2k)!}M_{2k}(E).
\]
\end{theorem}

\begin{proof}
We explain the proof only in the case when there exists a fiberwise oriented GMF $f_t$ on $E$. Then $-f_t$ is also a fiberwise oriented GMF. The negative eigenspace bundle of $-D^2f_t$ is equal to the positive eigenspace bundle $\eta$ for $D^2f_t$. Also, the index of the critical points change from $i$ to $n-i$. So, the framing principle gives:
\[
	\t_{2k}(E)=\t_{2k}(C(f_t))+(-1)^k\zeta(2k+1)p_\ast^\Sig(ch_{2k}(\xi))
\]
\[
	\t_{2k}(E)=\t_{2k}(C(-f_t))+(-1)^{n+k}\zeta(2k+1)p_\ast^\Sig(ch_{2k}(\eta))	
\]
However, the algebraic torsion has the property that $\t_{2k}(C(-f_t))=(-1)^n\t_{2k}(C(f_t))$. This is the \emph{involution property}. Therefore, for $n=\dim F$ even, we get the following by adding the above two equations.
\[
	2\t_{2k}(E)=(-1)^k\z(2k+1)p_\ast^\Sig(ch_{2k}(T^vE))=(-1)^k\z(2k+1)\left(
	\frac1{(2k)!}M_{2k}(E)
	\right)
\]
The theorem follows by dividing by 2.
\end{proof}

This theorem says that $M_{2k}$ is proportional to the even part of $\t_{2k}$. The following calculation shows that the odd part of $\t_{2k}$ has the same size (but opposite sign).

\begin{theorem}\cite{[I:BookOne]}\label{FR torsion of linear sphere bundles}
The $S^n$-bundle associated to the $SO(n+1)$-bundle $\xi$ over $B$
has higher FR-torsion
\[
    \t_{2k}(S^n(\xi))=(-1)^{n+k}\z(2k+1)ch_{4k}(\xi)
\]
\end{theorem}

\begin{proof}
It follows from the stability of $\t_{2k}$ that $\t_{2k}=0$ on all linear disk bundles over $B$. Since the $S^{n+1}$ bundle association to $\xi$ is the union of two disk bundles along their common boundary $S^n(\xi)$ we have by additivity (Theorem \ref{properties of FR torsion}.2) that
\[
	\t_{2k}(S^{n+1}(\xi))=2\t_{2k}(D^{n+1}(\xi))-\t_{2k}(S^n(\xi))=-\t_{2k}(S^n(\xi))
\] The left hand side is equal to $(-1)^{n+1+k}\z(2k+1) ch_{2k}(\xi)$ by Corollary \ref{cor:FR torsion for sphere bundle with section}). The theorem follows.
\end{proof}

\begin{cor}\label{cor:FR torsion satisfies Ic, IIc}
Higher Franz-Reidemeister torsion $\t_{2k}$ is a higher torsion invariant with $s_n=\frac12(-1)^{n+k}\z(2k+1)$.
\end{cor}

\begin{proof}
The additivity axiom follows from Theorem \ref{properties of FR torsion}.2. For the transfer axiom we need to prove the following for any smooth unipotent bundle $E\to B$ with closed fiber $F$ and any oriented linear sphere bundle $S^m(\xi)\to E$.
\begin{equation}\label{eq in proof of the transfer axiom for FR torsion}
    \t_{2k}(S^m(\xi))_B=\chi(S^m)\t_{2k}(E)+tr^E_B(\t_{2k}(S^m(\xi))_E)
\end{equation}
There are four cases depending of the parity of $m$ and $n=\dim F$.

If $n,m$ have the same parity then the transfer formula (\ref{eq in proof of the transfer axiom for FR torsion}) is equivalent to the transfer formula for $M_{2k}$ which we already proved. (If $n,m$ are both odd the RHS is zero for both $\t_{2k}$ and $M_{2k}$.)

If $n$ is even and $m$ is odd then, by the previous case we know that
\begin{equation}\label{one and a half-th equation}
    \t_{2k}(S^{m+1}(\xi))_B=2\t_{2k}(E)+tr^E_B(\t_{2k}(S^{m+1}(\xi))_E).
\end{equation}
But, the LHS of this is
\begin{equation}\label{second eq in proof of the transfer axiom for FR torsion}
    \t_{2k}(S^{m+1}(\xi))_B=2\t_{2k}(D^{m+1}(\xi))_B-\t_{2k}(S^m(\xi))_B
    = 2\t_{2k}(E)-\t_{2k}(S^m(\xi))_B
\end{equation}
by additivity. Combining (\ref{one and a half-th equation}) and (\ref{second eq in proof of the transfer axiom for FR torsion}) we get
\[
	\t_{2k}(S^m(\xi))_B=-tr^E_B(\t_{2k}(S^{m+1}(\xi))_E).
\]
However,
\[
    \t_{2k}(S^{m+1}(\xi))_E=-\t_{2k}(S^{m}(\xi))_E
\]
by Theorem \ref{FR torsion of linear sphere bundles}. So, (\ref{eq in proof of the transfer axiom for FR torsion}) holds in this case.

In the case $n$ odd and $m$ even, the LHS of (\ref{second eq in proof of the transfer axiom for FR torsion}) is zero since it is a multiple of
\[
	M_{2k}(S^{m+1}(\xi))_B=   \chi(S^{m+1})M_{2k}(E)+tr^E_B(M_{2k}(S^{m+1}(\xi))_E)=0+0.
\]
and the second term on the RHS of (\ref{eq in proof of the transfer axiom for FR torsion}) is zero since the transfer map is zero. So Equation (\ref{eq in proof of the transfer axiom for FR torsion}) holds proving the transfer axiom in all four parity cases.
\end{proof}

\begin{cor}
The odd part of higher FR-torsion is given by
\[
    \t_{2k}^-(E)=\t_{2k}(E)
    -(-1)^k\frac{\z(2k+1)}{2(2k)!}M_{2k}(E).
\]
\end{cor}

\begin{remark}A corollary of the main theorem is that any odd torsion invariant must be proportional to the above expression.
In particular it is expected that the nonequivariant analytic torsion classes of \cite{[Bismut-Lott95]} and \cite{[BG2]} are odd torsion invariants and
therefore proportional to $\t_{2k}^-$.
\end{remark}




\end{document}